\documentclass[a4paper,11pt]{article}

\title{Paralic confinement computations in coastal environment with interlocked areas}
\date{\today}
\author{
Jean-Philippe \textsc{Bernard}\thanks{Virtual Plants, C.C. 06002, 95 rue de la Gal\'era, 34095 Montpellier Cedex 5, France.},
Emmanuel \textsc{Fr\'enod}\thanks{Universit\'e Europ\'eenne de Bretagne, LMBA (UMR CNRS 6205), Universit\'e de Bretagne-Sud \& Inria Nancy-Grand Est, CALVI Project.}, 
Antoine \textsc{Rousseau}\thanks{Inria \& LJK, 95 rue de la Gal\'era, 34090 Montpellier, France. A. Rousseau is supported by Labex NUMEV and the Lefe-Insu project CoCoA.}
}

\usepackage[T1]{fontenc}		
\usepackage[hmargin=2.5cm,vmargin=3cm]{geometry}
\usepackage{fancyhdr}			
\usepackage{lastpage}				
\usepackage{graphicx}			
\usepackage{amsmath}			
\usepackage{amsfonts}			
\usepackage{amssymb}			
\usepackage{dsfont}				
\usepackage{stmaryrd}			
\usepackage{multirow}			
\usepackage{tabularx}				
\usepackage{listings}				
\usepackage{array}					
\usepackage{xcolor}				
\usepackage{ulem}					
\usepackage{hyperref}				
\usepackage{tikz}					
\usepackage{subfig}				
\usepackage{cases}	
\graphicspath{{img/}}

\definecolor{marron}{RGB}{116,54,0}

\newcommand{\be}{\begin{equation}}
\newcommand{\ee}{\end{equation}}
\newcommand{\bea}{\begin{eqnarray}}
\newcommand{\eea}{\end{eqnarray}}
\newcommand{\frad}[2]{\frac{\displaystyle {#1}}{\displaystyle {#2}}}

\newcommand{\dn}[1]{ \frad{\partial {#1}}{\partial n }}
\newcommand\1{\leavevmode\hbox{\rm \small1\kern-0.35em\normalsize1}}

\newcommand{\modifjp}[1]{#1}

\newtheorem{rmq}{Remark}

\begin{document}

\maketitle

\section{Abstract}

This paper is in the continuity of a work program, initiated in {Fr\'enod \& Goubert \cite{FreGou2007}, {Fr\'enod \& Rousseau \cite{FrRo12}} and 
{Bernard, Fr\'enod \& Rousseau \cite{bernard:hal-00776060}}. Its goal is to develop an approach of the paralic confinement usable from
the modeling slant, before implementing it in numerical tools.

\modifjp{More specifically, we here deal with the multiscale aspect of the confinement. If a paralic environment is separated into two (or more) connected areas, we will show that is possible to split the confinement problem into two related problems, one for each area.
At the end of this paper, we will focus on the importance of the interface length between the two subdomains.}

\section{Introduction}

Paralic confinement is one of the most pertinent parameters controlling the features of living species in paralic environments  ({\it i.e.} environments such as lagoons, estuaries, bays, {\it etc.}).
It was first introduced by {Gu\'elorget \& Perthuisot \cite{GP,GuPe83}} and it is linked with  nutrient concentration of water in the paralic environments.
It was widely discussed and tested in {Gu\'elorget, Frisoni \&  Perthuisot \cite{GuelEtAl}},
{Gu\'elorget \textit{et al.} \cite{GuGaLoPe}},
{Ibrahim \textit{et al.} \cite{IbGuFrRoMaPe}},
{Debenay, Perthuisot \& Colleuil \cite{DePerColl1993}},
{Redois \& Debenay \cite{RedDeb1996}},
{Barnes \cite{Barnes1994}},
 {Fr\'enod \& Goubert \cite{FreGou2007}}
 and {Tagliapietra \textit{et al.} \cite{Traga2009}}.\\
 Its knowledge in a given paralic environment is an important factor for supporting decision of decision makers acting on the paralic environment.
 For instance, it can be used to help and choose the exact localization of shellfish farms, or to estimate the impact of the building of a new dyke or dam.
 \\
 
 Since the recent works of {Fr\'enod \& Goubert \cite{FreGou2007}, {Fr\'enod \& Rousseau \cite{FrRo12}} and  {Bernard, Fr\'enod \& Rousseau \cite{bernard:hal-00776060}},
 we know that it is possible to develop a methodology to simulate numerically the confinement in any paralic environment.
We now enter  a phase of our work program which long term objective is to provide an operational tool to compute the paralic confinement in any point of any 
 paralic environment on earth, only from bathymetry and oceanographic data.
 \\
 Many questions need to be reached before achieving such an objective and we start here by tackling a modest (but important) one, which is related to the capability
 of computing separately paralic confinement in two connected areas of a given paralic environment. We shall particularly focus on the interface boundary condition that is required for such a coupling. 

\section{Modeling issues}\label{sec:model}
\subsection{Interlocked areas} \label{subset:explain} 
Coastal environment is made of interlocked areas. For instance, if we take a look at a marsh in a Mediterranean lagoon,
we face with the following cascade of areas:  Atlantic Ocean~- Gibraltar Strait~- Mediterranean sea - lagoon entrance - the lagoon - marsh entrance - the marsh.

Beside this, coastal environments may present a wide range of scales. For instance in the previously evoked cascade,
the marsh is several tens of meters large, while the lagoon size is about ten kilometers. Those two scales are small when compared with
the characteristic size of the Mediterranean sea, which is itself small with respect to the Atlantic Ocean dimension. 
Sizes of transition inlets -~Gibraltar Strait  and the lagoon entrance~- need also to be taken into account.\\

As it will be recalled in the sequel (see also \cite{FrRo12}),
the numerical computation of the paralic confinement in coastal environment first requires  the computation of the water flow essentially going from the ocean to the
coastal environment far end, induced by the combined effect of evaporation, tide and fresh water inputs from the rivers. Once the flow is known, a tracer
following this flow and undergoing diffusion, is then computed. At the end of the process, this tracer provides the value of the paralic confinement.

Because of the wide range of scales appearing in coastal environments (see above), a confinement simulation may rely on several mathematical models that one needs to couple. In the case where the two coupled models are identical, we face a domain decomposition problem.\\
Even if the coupling between the two (ore more) subdomains should actually be two-way, we will focus on one-way exchanges (from the open deep sea to the lagoon, and finally to the marsh). In other words, when we decompose a computational domain in two parts, we will consider a main lagoon (the part that is directly connected to the open sea) and a secondary lagoon (see Figure \ref{fig:domaines} below). The confinement field in the main lagoon will have to be computed accurately in a truncated domain $\Omega^{\textrm{main}}=\Omega\,\backslash\,\Omega^{\textrm{seg}}$, while the simulation in the secondary lagoon 
$\Omega^{\textrm{seg}}$ is nothing but another ``classical'' (\textit{i.e.} monodomain) confinement simulation, with the main lagoon playing the role of the open deep sea.\\
In the next subsection we recall the model previously developed to compute the paralic confinement field
in a lagoon. In subsection \ref{sec:multiscale}, we introduce the domain decomposition $\Omega=\Omega^{\textrm{main}}\cup\Omega^{\textrm{seg}}$ and the related issues. We pay a particular attention to the interface (and the related boundary conditions) between the main and secondary lagoons.

\begin{figure}[!h]
	\centering 
\begin{tikzpicture}[scale=1.2]
\def\RM{1}
\def\RS{0.6}
\def\CL{0.5}
\def\CPrct{50}
\def\EPrct{75}

\def\CW{\CPrct/100*\RM};
\def\EW{\EPrct/100*\RS};

\def\ACM{asin(\CW/\RM/2)}
\def\ACS{asin(\CW/\RS/2)}
\def\AEM{asin(\EW/\RM/2)}

\coordinate (Eleft) at ({-\EW/2},{-\CL/2 - \RM*cos(\ACM) - \RM*cos(\AEM)});
\coordinate (Cleft) at ({-\CW/2},0);
\node (OM) at (0,{-\CL/2 - \RM*cos(\ACM)}) {$\Omega^\textrm{main}$};
\node (OS) at (0,{\CL/2 + \RS*cos(\ACS)}) {$\Omega^\textrm{seg}$};
\node (O) at ({-2*\RM},0) {$\Omega$};

\begin{scope}[thick]
	\draw (Eleft) arc ({270-\AEM:90+\ACM:\RM}) --++ (0,{\CL}) arc ({270-\ACS:-90+\ACS:\RS}) --++(0,{-\CL}) arc ({90-\ACM:-90+\AEM:\RM});
	\draw[green!66!black] (Eleft) --++ (\EW,0) node[midway,below]{$\Gamma^{in}$};
\end{scope}

\draw[red] (Cleft) --++ (\CW,0) node[right]{$\Gamma$};

\begin{scope}[blue]
	\draw[densely dotted] ({-\CW/2},{-\CL/2}) --++ (0,{-\RM*0.5}) ++ (\CW,0) -- ({\CW/2},{-\CL/2});
	\draw[latex-latex] ({-\CW/2},{-\CL/2-\RM*0.5}) --++ (\CW,0) node[midway,above]{$\delta$};
\end{scope}

\begin{scope}
	\draw (O) -- (OM);
	\draw (O) -- (OS);
	\draw (OM) ++ ({-\RM},0) node[left] {$\Gamma^0$};
\end{scope}
\end{tikzpicture}
\hspace*{2cm}
\begin{tikzpicture}[scale=1.2]
\def\RM{1}
\def\RS{0.6}
\def\CL{0.5}
\def\CPrct{25}
\def\EPrct{75}

\def\CW{\CPrct/100*\RM};
\def\EW{\EPrct/100*\RS};

\def\ACM{asin(\CW/\RM/2)}
\def\ACS{asin(\CW/\RS/2)}
\def\AEM{asin(\EW/\RM/2)}

\coordinate (Eleft) at ({-\EW/2},{-\CL/2 - \RM*cos(\ACM) - \RM*cos(\AEM)});
\coordinate (Cleft) at ({-\CW/2},0);
\node (OM) at (0,{-\CL/2 - \RM*cos(\ACM)}) {$\Omega^\textrm{main}$};
\node (OS) at (0,{\CL/2 + \RS*cos(\ACS)}) {$\Omega^\textrm{seg}$};
\node (O) at ({-2*\RM},0) {$\Omega$};

\begin{scope}[thick]
	\draw (Eleft) arc ({270-\AEM:90+\ACM:\RM}) --++ (0,{\CL}) arc ({270-\ACS:-90+\ACS:\RS}) --++(0,{-\CL}) arc ({90-\ACM:-90+\AEM:\RM});
	\draw[green!66!black] (Eleft) --++ (\EW,0) node[midway,below]{$\Gamma^{in}$};
\end{scope}

\draw[red] (Cleft) --++ (\CW,0) node[right]{$\Gamma$};

\begin{scope}[blue]
	\draw[densely dotted] ({-\CW/2},{-\CL/2}) --++ (0,{-\RM*0.5}) ++ (\CW,0) -- ({\CW/2},{-\CL/2});
	\draw[latex-latex] ({-\CW/2},{-\CL/2-\RM*0.5}) --++ (\CW,0) node[midway,above]{$\delta$};
\end{scope}

\begin{scope}
	\draw (O) -- (OM);
	\draw (O) -- (OS);
	\draw (OM) ++ ({-\RM},0) node[left] {$\Gamma^0$};
\end{scope}
\end{tikzpicture}
	\caption{\label{fig:domaines}Computational domain $\Omega$. The lower part (below the interface $\Gamma$) is denoted $\Omega^\textrm{main}$. The upper part is denoted $\Omega^\textrm{seg}$. \underline{Left}: large interface ($\delta=20r_{1}/100$), where $r_{1}$ is the radius of $\Omega^\textrm{main}$. \underline{Right}: tinier interface ($\delta=10r_{1}/100$).\\
	Simulations will be conducted with $\delta/r_{1}=5\%, 10\%, 15\%$ and $20\%$.}
\end{figure}
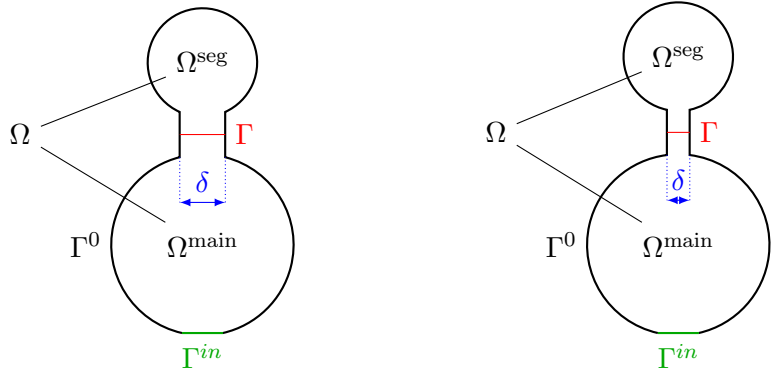

\subsection{Domain and equations for paralic confinement computation in a lagoon} \label{subsec:confmodels} 
We consider (see Figure \ref{fig:Schema}) a lagoon that  is a cylinder with base a regular, connected and bounded domain $\Omega\subset\mathbb{R}^2$ with boundary $\partial\Omega$. This boundary is shared into $\Gamma^\text{in}$ and $\Gamma^{0}$ with $\Gamma^\text{in}\cap\Gamma^{0}=\emptyset$. Any point in $\overline\Omega$ is denoted $(x,y)$. The lagoon seabed is described by a piecewise continuous function $b:\Omega\longrightarrow\mathbb{R}^+$, where $b(x,y)$ represents the bathymetry level at the horizontal position $(x,y)\in\Omega$. The water altitude $h$ is such that $h>\sup_{\Omega}\{b\} $, exluding outcrops. In summary, the geometrical model of the lagoon writes:
\begin{equation}
\label{lagoon}
\mathcal{L}agoon=\Big\{(x,y,z), \; (x,y)\in\Omega, \; b(x,y)<z<h\Big\}.
\end{equation}

\begin{figure}[!h]
	\centering 
	\includegraphics[width=8.2cm]{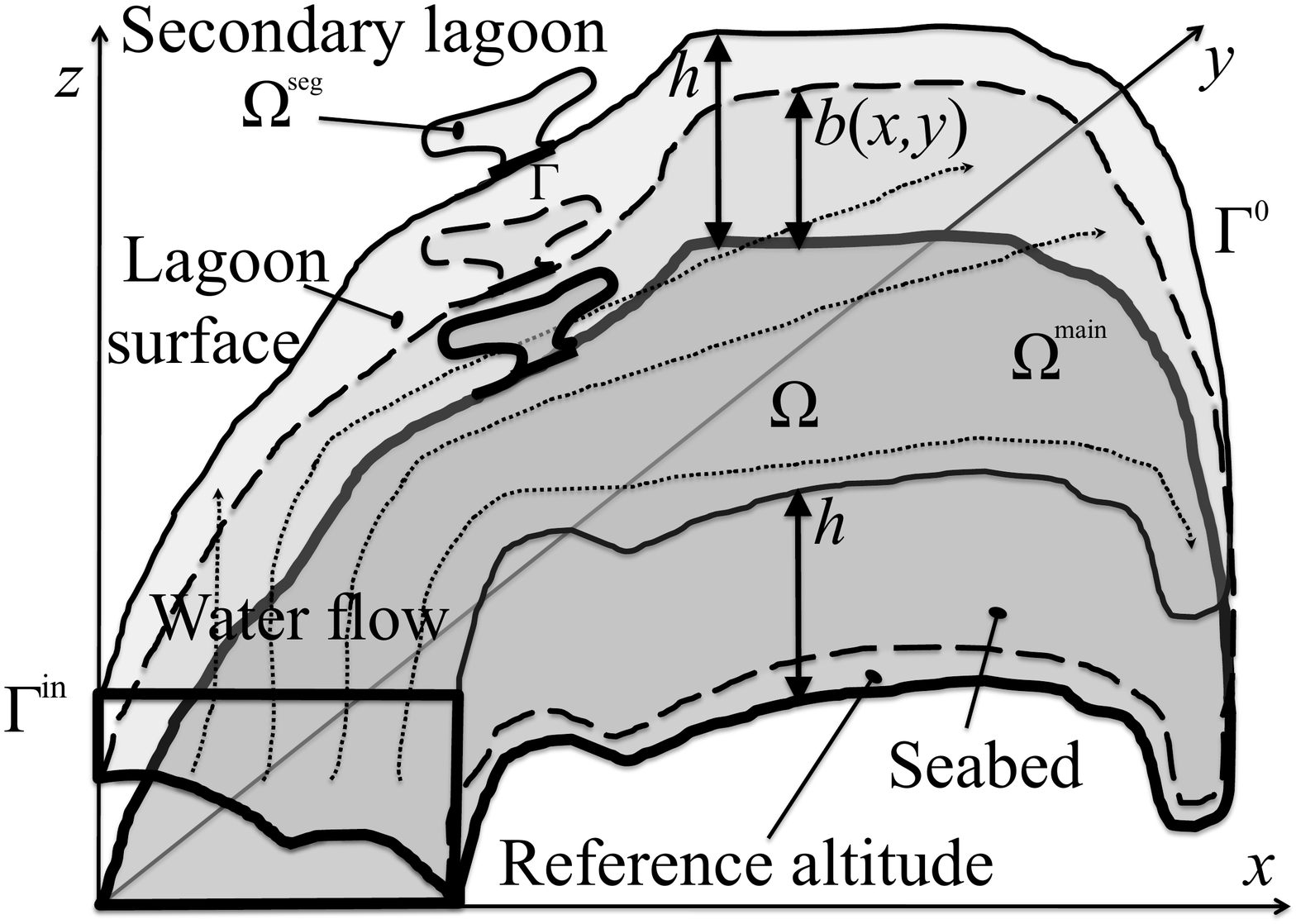}
	 \includegraphics[width=6.9cm]{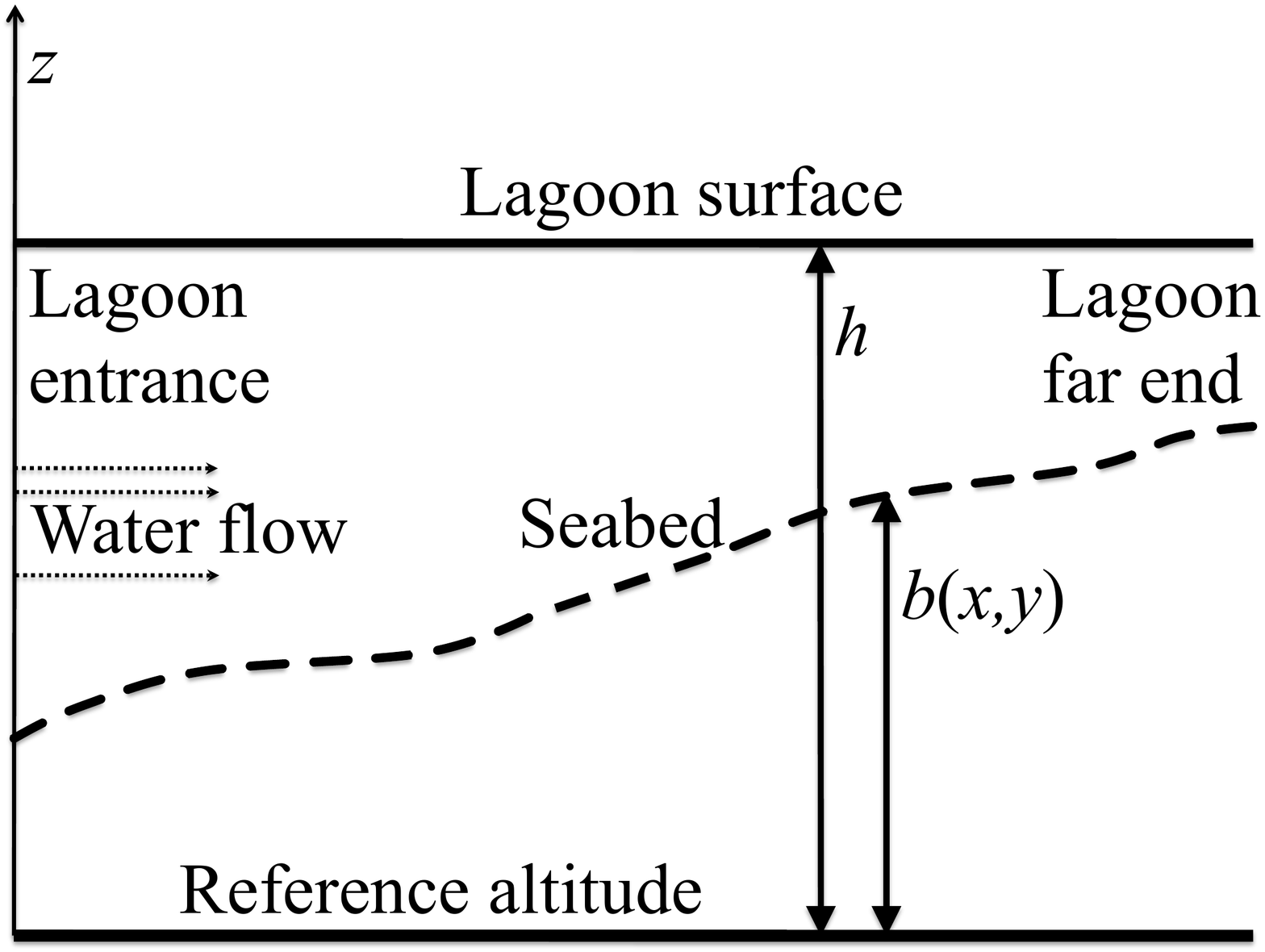}\vspace{-5mm}
	\caption{\label{fig:Schema}Left: lagoon geometry,  including a secondary lagoon. Right: A section of the lagoon geometry over a line going from the lagoon entrance to
	                                                the lagoon far end.}
\end{figure}
The  instantaneous confinement 
\modifjp{is - at a given time $t$ and a given position $(x,y)$ - the amount of time the water particle located at the position $(x,y)$ at time $t$ has spent inside the lagoon water mass.}
It is related to nutrient concentration at position $(x,y)$ and time $t$. 
It is the result of two phenomena. In the first place, evaporation process 
generates a flow from the ocean to the lagoon far end.
Essentially, ocean waters can be seen as having a very high
nutrient concentration and when those waters travel towards the lagoon, they meet living species that take off nutrients. Then along the travel, their nutrient concentration decreases.
Secondly, the nutrient concentration undergoes diffusion because of water molecular  and eddy viscosities, but more pregnantly,  because of the chemical precesses involved in 
the dissolution phenomena. This diffusivity is small, when compared with the average consequence of the flow. Nevertheless, it is important in places where the velocity of the flow is
small, especially in the lagoon far end.
\\
Consequently, in order to compute the instantaneous confinement, we use a passive tracer $g_{t}$ advected by the water velocity field $\pmb{u}$ and undergoing diffusion. As shown in {Fr\'enod \& Rousseau \cite{FrRo12}}, this model  is compatible with any lagoon geometry (shape and bathymetry), with the only restriction that intertidal zones and seabed outcrops are not taken into account. 
The idea developped in \modifjp{\cite{bernard:hal-00776060} is to solve the following advection-diffusion problem:}
for any time $t>0$ and given a sufficiently large time $T$, the solution $g_t=g_t(\tau,x,y)$ of
\begin{equation}\label{EqConfVisc}
\modifjp{\left\{
\begin{array}{lrcl}
\lefteqn{\forall 0<\tau<T, \forall (x,y)\in\Omega,}&&&\\
&\dfrac{\partial g_t}{\partial \tau}(\tau,x,y)+\pmb{u}(t-T+\tau,x,y)\cdot\nabla g_t(\tau,x,y) - \nu\Delta g_t(\tau,x,y) &=& 0,\\[0.35cm]
\lefteqn{\forall 0<\tau<T,\forall (x,y)\in{\Gamma}^\text{in},}&&&\\
&g_t(\tau,x,y) &=& T-\tau, \\[0.2cm]
\lefteqn{\forall (x,y)\in\Omega,}&&&\\
&g_t(0,x,y) &=& T,
\end{array}
\right.}
\end{equation}
is such that $g_t(T,x,y)$ is a good approximation of the value of the instantaneous confinement at time $t\in\mathbb R_+$ and position $(x,y)\in\Omega$.
\modifjp{Here $\tau$ is a variable related to the time to spend into the lagoon and} $\nu$ is the small (when compared with the average value of $\pmb{u}$) diffusivity coefficient
\modifjp{that models the nutrient hability to spread out the water.}
 In system  
\eqref{EqConfVisc}, the water velocity field $\pmb{u}(t,x,y)$ 
\modifjp{may be induced by several phenomena (such as evaporation, tide, river input, etc.), which are modeled by the generic function $\theta$. We may compute this field}
by solving the following equation:
\begin{equation}\label{eq:velocity} 
\left\{
\begin{array}{rcll}
-\nabla\cdot\big[(h-b)\pmb{u}\big](t,x,y) &=& \theta(t,x,y), & \;\forall t>0, \forall (x,y)\in\Omega,\\
\nabla\times\pmb{u}  &=&  0 , & \;\forall t>0, \forall (x,y)\in\Omega,\\
\pmb{u}\cdot{n} &=& F^{\text{in}}(t,x,y), & \;\forall t>0, \forall (x,y)\in\Gamma^\text{in},\\
\pmb{u}\cdot{n} &=& 0, & \;\forall t>0, \forall (x,y)\in\Gamma^{0},
\end{array}
\right.
\end{equation}
where $n$ stands for the unitary vector orthogonal to $\partial\Omega$  pointing outside $\Omega$ and where $F^{\text{in}}$ is a function defined on ${\Gamma}^\text{in}$ such that :
\begin{equation}
\label{defFForPlain} 
\int_{\modifjp{\Gamma^{\text{in}}}}\big[(h-b)F^{\text{in}}\big](t,x,y)\:\text dl = \int_\Omega\theta(t,x,y)\:\text dx\:\text dy, \;\forall t>0.
\end{equation}

\modifjp{We remark that} the velocity field $\pmb{u}$ can be separated in several ``elementary'' velocity fields, each of those being solely induced by one single process
Consequently, depending on those processes, the function $\theta$ can model
\modifjp{one or}
several phenomena.
The sytem of equations \eqref{eq:velocity} is solved thanks to its corresponding velocity potential formulation (see \cite{FrRo12}). Provided that $\nabla\times u=0$, we write $u=\nabla\psi$ and solve the following Laplace equation for $\psi$:
\begin{subnumcases}{\label{eq:psi}}
-\nabla\cdot\big[(h-b)\nabla\psi\big]=\theta\quad\mbox{on } \Omega,\\
\dn{\psi}=F^{\text{in}}\quad\mbox{ on } \Gamma^\text{in},\\
\dn{\psi}=0\quad\mbox{ on } \Gamma^ {0}.
\end{subnumcases}
In subsection \ref{sec:multiscale} and section \ref{sec:numerics} below we will consider equations for the velocity in truncated regions of the lagoon (\textit{e.g.} equations \eqref{eq:velocitySUB} and \eqref{eq:velocityMAIN}). Naturally, these equations will always be solved thanks to the potential formulation, even if this is not explicitly specified.
\subsection{Domain decomposition for a lagoon with a secondary lagoon}\label{sec:multiscale} 
When in a lagoon a clearly separated entity -~so called secondary lagoon~- exists, we want to \modifjp{split} problem  \eqref{EqConfVisc},   \eqref{eq:velocity} and \eqref{defFForPlain}
into two problems -~a first one set in the secondary lagoon and another one set in the remainder of the lagoon~- being connected by conditions on the secondary lagoon entrance.
Naturally, we want the concatenation of the results to approximate a solution of the system  \eqref{EqConfVisc},   \eqref{eq:velocity} and \eqref{defFForPlain}
set in the whole lagoon with a good accuracy.

The way to account for this situation, in what concerns the geometrical aspects, consists in sharing
lagoon $\Omega$ into three parts (see the left picture in Figure \ref{fig:Schema}): the secondary lagoon $\Omega^\textrm{seg}$, 
the main part  of the lagoon $\Omega^\textrm{main}$ and their common boundary $\Gamma$. They are
such that $\Omega^\textrm{seg}$ and  $\Omega^\textrm{main}$ are open subsets of $\Omega$,  
$\Gamma = \partial\Omega^\textrm{seg}\cap\partial\Omega^\textrm{main}$ and
$\Omega^\textrm{seg}\cup \Omega^\textrm{main}\cup\Gamma =\Omega$.
Moreover we assume that the secondary lagoon is not near the lagoon entrance, which is translated by
$\partial\Omega^\textrm{seg}\cap\modifjp{\Gamma^{\text{in}}}=\emptyset$.\\
~\\

The key problem -~and the most difficult one~- is the obtention of the water flow and of the tracer within the main part
of the lagoon without computing them in the secondary lagoon. We now focus on this issue and, in a first place, we build the system of equations to compute the velocity of the water flow in the 
main part of the lagoon $\Omega^\textrm{main}$.
The two first equations of problem \eqref{eq:velocity} are retained, because they describe the physics of water transport. The third equality of \eqref{eq:velocity} that translates that water cannot escape from the secondary lagoon through the 
shore will also be kept. Thirdly,  a condition on the interface -~that we will discuss hereafter~- will be written.
Hence, in the main part of the lagoon $\Omega^\textrm{main}$, we write the following problem that allows us to obtain the velocity
\begin{equation}\label{eq:velocityMAIN} 
\left\{
\begin{array}{rcll}
-\nabla\cdot\big[(h-b)\pmb{u}^\textrm{main}\big](t,x,y) &=& \theta(t,x,y), & \;\forall t\in\mathbb R, \forall (x,y)\in\Omega^\textrm{main},\\
\nabla\times\pmb{u}^\textrm{main}  &=&  0 , & \;\forall t\in\mathbb R, \forall (x,y)\in\Omega^\textrm{main},\\
\pmb{u}^\textrm{main}\cdot{n} &=& F^{\text{in}}(t,x,y), & \;\forall t\in\mathbb R,\forall (x,y)\in{\Gamma}^\text{in},\\
\pmb{u}^\textrm{main}\cdot{n} &=& 0, & \;\forall t\in\mathbb R,\forall (x,y)\in\Gamma^0,\\
\pmb{u}^\textrm{main}\cdot{(-n^\text{trans})} &=& - F(t,x,y), & \;\forall t\in\mathbb R,\forall (x,y)\in\Gamma,\\
\end{array}
\right.
\end{equation}
where vector $n$ and function $F^{\text{in}}$  have the same definitions as in system \eqref{eq:velocity}  and equality \eqref{defFForPlain}, 
where $n^\text{trans}$ stands for the unitary vector, orthogonal to $\Gamma$, and pointing inside $\Omega^\textrm{main}$ 
(or outside $\Omega^\textrm{seg}$). The function $F$ is to be determined.\\
It is clear that if the solution $\pmb{u}$ of \eqref{eq:velocity} were known, we would choose 
$F= \pmb{u}\cdot{n^\text{trans}}$ and then we would obtain a solution $\pmb{u}^\textrm{main}$
of  \eqref{eq:velocitySUB} that would be such that $\pmb{u}^\textrm{main} = {\pmb{u}}_{|\Omega^\textrm{main}}$.
Yet, we work under the assumption that the solution of  \eqref{eq:velocitySUB} is not known (we indeed want to compute
$\pmb{u}^\textrm{main}$ to have the value of ${\pmb{u}}$ in ${\Omega^\textrm{main}}$).  Anyway,  using the Laplace-Neumann
compatibility condition (the same that brings us to write \eqref{defFForPlain} for system  \eqref{eq:velocity} for $\pmb{u}$), we know
an information on $F$ which is 
\begin{equation}
\label{CompForFTrans} 
\int_{\Gamma}\big[(h-b)F\big](t,x,y)\:\text dl = \int_{\Omega^\textrm{seg}}\theta(t,x,y)\:\text dx\:\text dy, \;\forall t\in\mathbb R,
\end{equation}
and that translates that the quantity of water entering $\Omega^\textrm{seg}$ -~and so leaving  $\Omega^\textrm{main}$~- through $\Gamma$ compensates for what is consumed
by the process modeled by $\theta$ over $\Omega^\textrm{seg}$.
Knowing this information, we can consider that the missing information is the profile (or the shape) of $F$ along the interface $\Gamma$.
Approximations of this profile, that are classical and known as giving proper results,  can be used. For instance, we can chose $F$ as being constant 
along $\Gamma$ or being 
a Poiseuille profile (see \cite{poiseuille1844mouvement}).\\

Having $\pmb{u}^\textrm{main}$ on hand, and then considering that we consequently know $\pmb{u}$ on $\Omega^\textrm{main}$
with a good accuracy, in order to compute the passive tracer given by  \eqref{EqConfVisc} only on $\Omega^\textrm{main}$, we will consider 
the following problem:
\begin{equation}\label{EqConfViscMAIN} 
\modifjp{\left\{
\begin{array}{lrcl}
\lefteqn{\forall 0<\tau<T, \forall (x,y)\in\Omega,}&&&\\
&\dfrac{\partial g_t^\textrm{main}}{\partial \tau}(\tau,x,y)+\pmb{u}^\textrm{main}(t-T+\tau,x,y)\cdot\nabla g_t^\textrm{main}(\tau,x,y) - \nu\Delta g_t^\textrm{main}(\tau,x,y) &=& 0,\\[0.35cm]
\lefteqn{\forall 0<\tau<T,\forall (x,y)\in{\Gamma}^\text{in},}&&&\\
&g_t^\textrm{main}(\tau,x,y) &=& T-\tau, \\[0.2cm]
\lefteqn{\forall (x,y)\in\Omega,}&&&\\
&g_t^\textrm{main}(0,x,y) &=& T,
\end{array}
\right.}
\end{equation}
which straightforwardly comes from \eqref{EqConfVisc} replacing $\pmb{u}$ by $\pmb{u}^\textrm{main}$.
This system has to be coupled with boundary conditions on the interface $\Gamma$.
As the diffusivity coefficient $\nu$ is small, following Halpern  \cite{Ha86}, choosing
the following Neumann condition
\begin{gather}
\label{ABC0}
\frac{\partial g^\textrm{main}_t}{\partial (-n^\text{trans})}(\tau,x,y) = 0, ~~~ \forall 0<\tau<T, ~~\forall (x,y)\in\Gamma,
\end{gather}
will give a solution $g_t^\textrm{main}$ which will correctly approach $g_t$ over $\Omega^\textrm{main}$.\\
~\\

Once brought a way to tackle the key problem, we can notice that we can implement a way to obtain the water flow and 
the tracer within the secondary lagoon.
As a decision was made concerning the  profile of $F$, and so concerning $F$ on $\Gamma$,
we can write the following system to obtain the velocity field $\pmb{u}^\textrm{seg}$ in the secondary lagoon $\Omega^\textrm{seg}$:
\begin{equation}\label{eq:velocitySUB} 
\left\{
\begin{array}{rcll}
-\nabla\cdot\big[(h-b)\pmb{u}^\textrm{seg}\big](t,x,y) &=& \theta(t,x,y), & \;\forall t\in\mathbb R, \forall (x,y)\in\Omega^\textrm{seg},\\
\nabla\times\pmb{u}^\textrm{seg}  &=&  0 , & \;\forall t\in\mathbb R, \forall (x,y)\in\Omega^\textrm{seg},\\
\pmb{u}^\textrm{seg}\cdot{n} &=& 0, & \;\forall t\in\mathbb R,\forall (x,y)\in\Gamma^0,\\
\pmb{u}^\textrm{seg}\cdot{n^\text{trans}} &=& F(t,x,y), & \;\forall t\in\mathbb R,\forall (x,y)\in\Gamma,\\
\end{array}
\right.
\end{equation}
where $n$ has the same definition as in system  \eqref{eq:velocity} and where  $n^\text{trans}$ and $F$ are the one set 
to solve system \eqref{eq:velocityMAIN}.

On the other hand, $g_t^\textrm{main}$ on interface $\Gamma$ can be computed as a result of system  \eqref{EqConfViscMAIN} 
and gives the value of the tracer on this
interface. This function can be used as a Dirichlet boundary condition for the problem giving $g_t^\textrm{seg}$, which will be
close to $g_t$ in $\Omega^\textrm{seg}$. This problem reads:
\begin{equation}\label{EqConfViscSUB}
\modifjp{\left\{
\begin{array}{lrcl}
\lefteqn{\forall 0<\tau<T, \forall (x,y)\in\Omega,}&&&\\
&\dfrac{\partial g_t^\textrm{seg}}{\partial \tau}(\tau,x,y)+\pmb{u}^\textrm{seg}(t-T+\tau,x,y)\cdot\nabla g_t^\textrm{seg}(\tau,x,y) - \nu\Delta g_t^\textrm{seg}(\tau,x,y) &=& 0,\\[0.35cm]
\lefteqn{\forall 0<\tau<T,\forall (x,y)\in{\Gamma}^\text{in},}&&&\\
&g_t^\textrm{seg}(\tau,x,y) &=& T-\tau, \\[0.2cm]
\lefteqn{\forall (x,y)\in\Omega,}&&&\\
&g_t^\textrm{seg}(0,x,y) &=& T,
\end{array}
\right.}
\end{equation}
with
\begin{gather}
\label{BCforSUB}
g^\textrm{seg}_t(\tau,x,y) = g^\textrm{main}(\tau,x,y), ~~~ \forall 0<\tau<T, ~~\forall (x,y)\in\Gamma.
\end{gather}

\section{Numerical simulations}\label{sec:numerics}
\modifjp{In this section, we present numerical simulations of the lagoon described in Figure \ref{fig:domaines} (see Section \ref{sec:model}). These simulations were performed  with the finite element method implemented in the FreeFem++ software \cite{freefem++}. As in \cite{bernard:hal-00776060} the velocity equation is solved thanks to a Laplace equation on the velocity potential ($\pmb u = \nabla \psi$ and we use $P2$ elements for $\psi$), whereas the advection-diffusion equation on $g_t$ is solved thanks to $P1$ elements. }

\noindent We consider four different configurations, in order to enhance the importance of the interface width~$|\Gamma|$. The lagoon $\Omega$ is equally split in two parts $\Omega^\textrm{main}$ and $\Omega^\textrm{seg}$, which is the most general (unfavorable) case. Indeed, in cases where $|\Omega^\textrm{seg}|<|\Omega^\textrm{main}|$, the truncation error obtained in the numerical simulations is lower.\footnote{These simulations have been performed, but for the sake of clarity we only present here the case where $|\Omega^\textrm{seg}|=|\Omega^\textrm{main}|$.
}\\
For each of these configurations, we perform a numerical simulation of confinement in the whole domain $\Omega$, thanks to equations \eqref{EqConfVisc} and \eqref{eq:velocity}.
The corresponding numerical solutions will be considered as reference solutions and denoted $(\pmb u^\textrm{ref},g_{t}^\textrm{ref})$.

\subsection{Simulations without interface information}\label{sec:UC}
We now consider the numerical simulation of confinement in the truncated domain $\Omega^\textrm{main}$ (see Figure~\ref{fig:domaines}), in which we look for $\pmb u^\textrm{main},g_{t}^\textrm{main}$ solutions of systems \eqref{eq:velocityMAIN} and \eqref{EqConfViscMAIN}. For the sake of simplicity, we consider a flat bottom and set $h-b\equiv 1$, so that the boundary function $F$ is such that:

\begin{equation}\ 
F(t,x,y) = f_{pr}(t,x,y)\times\int_{\Omega^\textrm{seg}}\theta(t,x,y)\:\text dx\:\text dy \quad\forall t\in\mathbb R,\forall (x,y)\in\Gamma,\label{eq:BC1-u}
\end{equation}
where $f_{pr}$ denotes the (unknown) profile of the velocity along the interface $\Gamma$ and is such that $\int_{\Gamma}f_{pr}(t,x,y)d\sigma=1$. In the simulations below, we use a Poiseuille profile for $f_{pr}$.\\

\noindent We enumerate in Table \ref{tab:UC} the $L^{\infty}$ norm of the relative error
between the reference solution $g_{t}^\textrm{ref}$ restricted to the subdomain $\Omega^\textrm{main}$ and $g_{t}^\text{main}$ (computed in $\Omega^\textrm{main}$ from  \eqref{EqConfViscMAIN}-\eqref{ABC0}).
\modifjp{This error is defined by 
$$ \sup_{\Omega^{\text{main}}}\left|\frac{g_t^{\text{main}}-g_t^{\text{ref}}}{g_t^{\text{ref}}}\right|. $$}
\noindent We observe from the results in Table \ref{tab:UC} that the longer the interface, the larger the error. This is due to the lack of information we have at the interface (in particular on the velocity profile, see Section \ref{sec:model} and discussion above).

\begin{table}[!h]
\begin{minipage}{0.6\linewidth}
\begin{center}
\begin{tabular}{|c|c|c|c}
  \hline
  Configuration & Interface width $\delta/r_{1}$ & $L^{\infty}$ relative error\\
  \hline
  1 &  20\%& 0.0279855 \\
  2 &  15\%&  0.0212144\\
  3 &  10\%&  0.0133008\\
  4 &  5\%  &  0.00627107\\
  \hline
\end{tabular}
\end{center}
\end{minipage}
\begin{minipage}{0.3\linewidth}
\hspace*{1cm}
\includegraphics[width=\linewidth]{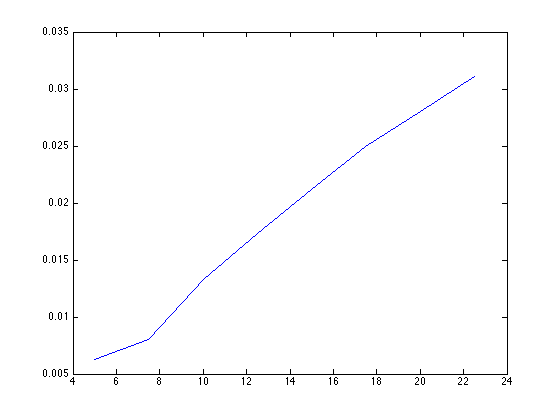}
\end{minipage}
\caption{\label{tab:UC}Relative error between confinements $g_{t}^\textrm{main}$ and $g_{t}^\textrm{ref}$ with the same diffusivity $\nu=0.01$ and four different interface widths. The error increases with the interface width. The figure only illustrates the results indicated in the table: the plot corresponds to the infinite relative error (column 3) as the function of the interface width (column 2).}
\end{table}
\subsection{Simulations with interface information}
We now reproduce the numerical simulations of Section \ref{sec:UC}, but instead of choosing a Poiseuille profile for the function $f_{pr}$, we use the exact profile provided by the knowledge of $\pmb{u}^\textrm{ref}$ on the interface\footnote{We could also imagine obtaining this information from measurements.}:
\begin{equation}
\pmb{u}^\textrm{main}\cdot{n^\text{trans}} = \pmb{u}^\textrm{ref}(t,x,y)\cdot{n^\text{trans}}\quad\forall t\in\mathbb R,\forall (x,y)\in\Gamma.
\end{equation}
\begin{rmq}
We obviously have that $\int_{\Gamma}\pmb{u}^\textrm{ref}(t,x,y)\cdot{n^\text{trans}}\;d\sigma=\int_{\Omega^\textrm{seg}}\theta(t,x,y)\:\text dx\:\text dy$, which means that the lack of information on $\Gamma$ in Equation \eqref{eq:BC1-u} concerns the velocity profile $f_{pr}$ rather than its average amplitude over the interface.
\end{rmq}
Then, thanks to the well-posedness of the corresponding velocity equations, uniqueness immediately insures that $\pmb u^\textrm{main}=\pmb u^\textrm{ref}_{|\Omega^\textrm{main}}$, that is to say the knowledge of the velocity profile along the interface $\Gamma$ insures the complete knowledge of the velocity in $\Omega^\textrm{main}$. As above, we enumerate the corresponding errors in Table \ref{tab:C}, and as expected the errors also depend (in the same manner) on the interface width, but are notably lower than those of Table \ref{tab:UC}.

\begin{table}[!h]
\begin{minipage}{0.6\linewidth}
\begin{center}
\begin{tabular}{|c|c|c|c}
  \hline
  Configuration & Interface width $\delta/r_{1}$ & $L^{\infty}$ relative error\\
  \hline
  1 & 20\% &  0.00267697\\
  2 & 15\% &  0.00204638\\
  3 &  10\%&  0.00134943\\
  4 &  5\%  &  0.000690316\\
  \hline
\end{tabular}
\end{center}
\end{minipage}
\hspace*{1cm}
\begin{minipage}{0.3\linewidth}
\includegraphics[width=\linewidth]{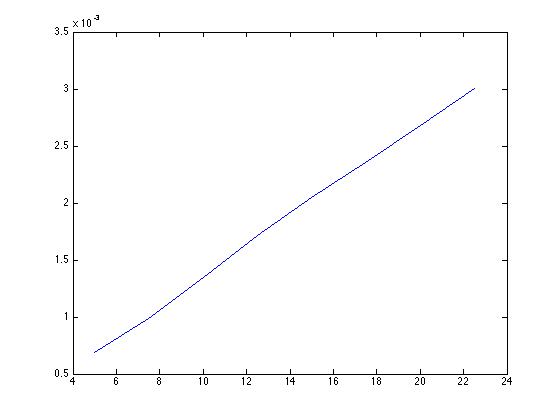}
\end{minipage}
\caption{\label{tab:C}Relative error between confinements $g_{t}^\textrm{main}$ and $g_{t}^\textrm{ref}$ with the same diffusivity $\nu=0.01$ and four different interface widths. The error increases with the interface width, and is notably lower than in Table \ref{tab:UC}. The figure only illustrates the results indicated in the table: the plot corresponds to the infinite relative error (column 3) as the function of the interface width (column 2).}
\end{table}

\noindent Furthermore, we can illustrate the quality of the Neumann boundary conditions \eqref{ABC0} used for the confinement equation with regard to the confinement diffusivity. Table \ref{tab:OBC-approx} illustrates that when the diffusivity is low, Neumann boundary conditions efficiently approximate the exact transparent boundary conditions, which was already proved in \cite{Ha86}.

\begin{table}[!h]
\begin{minipage}{0.6\linewidth}
	\centering 
	 \begin{tabular}{|c|c|}
	\hline
	Diffusivity 	$\nu$				& $L^{\infty}$ relative error \\\hline
	$1.10^{-1}$				& 0.010308 \\\hline
	$5.10^{-2}$	& 0.00323477 \\\hline
	$1.10^{-2}$				& 0.00267697 \\\hline
	$5.10^{-3}$	& 0.00195359 \\\hline
	\end{tabular}
\end{minipage}
\begin{minipage}{0.3\linewidth}
\includegraphics[width=\linewidth]{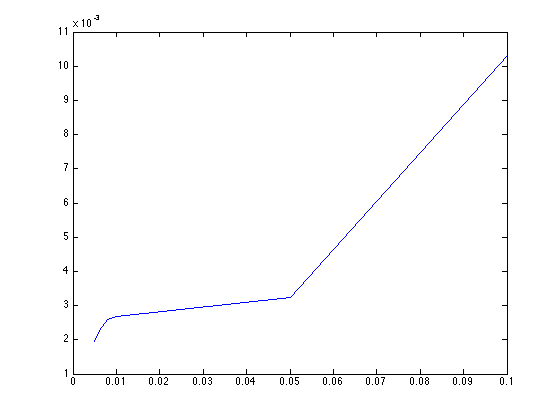}
\end{minipage}
\caption{\label{tab:OBC-approx}Relative error between confinements $g_{t}^\textrm{main}$ and $g_{t}^\textrm{ref}$ as a function of the diffusivity $\nu$, in configuration 4 (worst case, large interface width $\delta=20r_{1}/100$). The error increases with the diffusivity. The figure only illustrates the results indicated in the table: the plot corresponds to the infinite relative error (column 2) as the function of the diffusivity (column 1).}
\end{table}

\section{Conclusion}
In this paper we are interested in the truncation of computational domains in confinement models. This 
is a very important issue both for classical domain decomposition problems and for the numerical simulation of confinement in limited areas of large lagoons, which we consider here. As soon as the domain truncation is done, the most important question 
to address is the search for artificial boundary conditions  at the new boundary. It is known that their nature strongly depends on the PDE model that drives the considered process. Starting from the confinement model 
introduced in \cite{FrRo12}, we introduced some boundary conditions in order to limit the numerical error 
induced by the domain truncation.\\
The chosen confinement condition is a classical homogeneous Neumann boundary condition, which is known to be accurate for small diffusivity values (see \cite{Ha86}). The truncation error is actually mainly due to the lack of knowledge  of the velocity profile across the artificial boundary, particularly in the case where the interface and/or the secondary lagoon are large. It would be very interesting to evaluate how some partial informations on this velocity profile (provided by measurements) would improve the corresponding truncation error. We leave this to subsequent studies.

\section*{Acknowledgments}
The authors are very grateful to A. Fiandrino for fruitful discussions related to this article.
%
%

\bibliographystyle{alpha}
\bibliography{biblio}
\end{document}